\documentclass[11pt]{article}

\usepackage{amsmath}
\usepackage{amsthm}
\usepackage{amssymb}
\usepackage{mathptmx}
\usepackage{concrete}

\def\bel#1{\begin{equation}\label #1}

\def\l{\lambda}

\def\d{\delta}
\def\D{\Delta}

\def\G{\Gamma}

\def\nd{\noindent}
\def\p{\partial}

\def\BB{{\Bbb B}}

\def\<{\leq}
\def\>{\geq}

\begin{document}

\title{\bf Gradient estimate of an eigenfunction on a compact Riemannian manifold without boundary}
\author{Yiqian Shi $^{*\dagger}$ and Bin Xu$^{*\star}$}
\date{}
\maketitle

\nd{\small {\bf Abstract.}\ \ Let $e_\l(x)$ be an eigenfunction
with  respect to the Laplace-Beltrami operator $\Delta_M$ on a
compact Riemannian manifold $M$ without boundary: $\Delta_M
e_\l=\l^2 e_\l$. We show the following gradient estimate of
$e_\l$: for every $\l\geq 1$, there holds $\l\|e_\l\|_\infty/C\leq
\|\nabla e_\l\|_\infty\leq C\l \|e_\l\|_\infty$, where $C$ is a
positive constant depending only on $M$.  }

\footnote{$^*$ Department of Mathematics, University of Science
and Technology of China, Hefei 230026 China.\\
$^\dagger$ Supported in part by the National Natural Science
Foundation of China (No. 10671096).  e-mail: yqshi@ustc.edu.cn\\
$^\star$ The correspondent author is supported in part by the
National Natural Science Foundation of China (No. 10601053 and No.
10871184). e-mail: bxu@ustc.edu.cn}

\nd {\small {\bf Mathematics Subject Classification (2000):}\
Primary 35P20; Secondary 35J05}

\nd {\small {\bf Key Words:}\ Laplace-Beltrami operator,
eigenfunction, gradient estimate}

\section{Introduction}
\label{sec:sob1} \quad  Let $(M,\ g)$ be an $n$-dimensional
compact smooth Riemannian manifold without boundary and $\Delta_M$
the positive Laplace-Beltrami operator on $M$. Let $L^2(M)$ be the
space of square integrable functions on $M$ with respect to the
Riemannian density $dv(M)=\sqrt{{\bf g}(x)}\,dx:= \sqrt{\det\
(g_{ij})}\ dx$. Let $e_1(x),\,e_2(x),\,\cdots$ be a complete
orthonormal basis in $L^2(M)$ for the eigenfunctions of $\D_M$
such that $0= \l_0^2< \l_1^2\leq  \l_2^2\leq\,\cdots$ for the
corresponding eigenvalues, where $e_j(x)$ ($j=1,2,\cdots$) are
real valued smooth function on $M$ and $\l_j$ are nonnegative real
numers. Also, let ${\bf e}_j$ denote the projection onto the
1-dimensional space ${\bf C}e_j$. Thus , an $L^2$ function $f$ can
be written as $f=\sum_{j=0}^\infty {\bf e}_j(f)$, where the
partial sum converges in the $L^2$ norm. Let $\l$ be a positive
real number $\geq  1$. We define the spectral function $e(x,y,\l)$
and the unit band spectral projection operator $\chi_\l$ as
follows:
\[ e(x,y,\l):=\sum_{\l_j\leq \l} e_j(x) e_j(y)\ , \]
\[
{\chi}_\l f:=\sum_{\l_j\in (\l,\,\l+1]}{\bf e}_j(f)\ .\]

In 1968, H{\" o}rmander \cite{H1} obtained a one-term asymptotic
expansion of the spectral function of a positive definite elliptic
linear operator, whose Laplacian case, also called by the local
Weyl law, has the expression as follows:
\[ e(x,x,\l)=\l^n/\bigl((4\pi)^{n/2}\G(1+n/2)\bigr) +{\rm
O}(\l^{n-1}),\ \l\to\infty.
\]
As a consequence H{\" o}rmander proved the uniform estimate of
eigenfunctions for all $x\in M$: \[\sum_{\l_j\in (\l,\,\l+1]}
|e_j(x)|^2 \leq C\,\l^{n-1}\ .
\]
We note here that in the whole of this paper $C$ denotes a
positive constant which depends only on $M$ and may take different
values at different places. C. D. Sogge noted in (1.7) of
\cite{S4} the equivalence of the above estimate with the following
$L^\infty$ estimate of $\chi_\l$,
\begin{equation}
\label{equ:sup}
 ||{ \chi}_\l f||_\infty\leq C\l^{(n-1)/2} ||f||_2\
.
\end{equation}
Here $||f||_r$ ($2\leq r\leq \infty$) means the $L^r$ norm of the
function $f$ on $M$. As noted in Lemma 2.7 by the second author
\cite{XuB}, there also holds the similar equivalent relationship
between the following two gradient estimates for eigenfunctions
and $\chi_\l$:
\begin{equation}
\label{equ:spegrad} \sum_{\l_j\in (\l,\,\l+1]} |\nabla e_j(x)|^2
\leq C\,\l^{n+1}\quad \text{for all}\ x\in M,
\end{equation}
and
\begin{equation}
\label{equ:2infgrad} ||\nabla\,  { \chi}_\l f||_\infty\leq
C\l^{(n+1)/2} ||f||_2\ .
\end{equation}
Here $\nabla$ is the Levi-Civita connection on $M$.  In
particular, $\nabla f=\sum_j\, g^{ij}\partial f/\partial x_j$ is
the gradient vector field of a $C^1$ function $f$, the square of
whose length equals $\sum_{i,j} g^{ij}(\p f/\p x_i)(\p f/\p x_j)$,
where $(g^{ij})$ is the inverse of the metric matrix $g_{ij}$. By
using the wave group,  Yu Safarov and D. Vassiliev  proved a very
general theorem (Theorems 1.8.5 and 1.8.5 in \cite{SV}) on the
spectral function of a positive definite elliptic linear operator
so that the gradient estimate (\ref{equ:spegrad}) is its immediate
corollary. By using the parametix of the wave operator $\p^2/\p
t^2+\D_M$, the second author \cite{XuB} also proved in a slightly
different way the Laplacian case of the aforementioned theorems of
Yu Safarov and D. Vassiliev, which is also sufficient to deduce
the estimate (\ref{equ:spegrad}). X. Xu \cite{XuX} applied the
maximum principle argument to proving (\ref{equ:2infgrad}), and
used this estimate to give a new proof for the H{\" o}rmander
multiplier theorem. A. Seeger and C. D. Sogge \cite{SS} firstly
proved the the H{\" o}rmander multiplier theorem by using the
parametrix of the wave kernel. In this paper we stick at
estimating the gradient
of a single eigenfunction.\\

\nd {\bf Theorem} {\it Let $e_\l(x)$ be an eigenfunction with
respect to the positive Laplace-Beltrami operator $\Delta_M$ on an
$n$-dimensional compact smooth Riemannian manifold $M$ without
boundary: $\Delta_M e_\l=\l^2 e_\l$. Then, for every $\l\geq 1$,
there holds
\begin{equation}
\label{equ:grad} \l\|e_\l\|_\infty/C\leq \|\nabla
e_\l\|_\infty\leq C\l \|e_\l\|_\infty,
\end{equation}
where $C$ is a positive constant depending only on $M$.}\\

We can also obtain the above estimates (\ref{equ:spegrad}) and
(\ref{equ:2infgrad}) by H{\" o}rmander's local Weyl law and the
same argument of Theorem.  The details will be given in Section 4.
The second author announced a more general estimate for
$\|\nabla^k e_\l\|_\infty$ with $k$ a positve integers in Theorem
5.1 of \cite{XuB2}. However, the method we use in this paper can
not be applied to deducing higher derivative estimates. We plan to
discuss this question in a future paper.

We conclude the introduction by explaining the organization of
this paper. Sections 2-4 contain the proof of Theorem. In Section
2, we show the lower bound of the gradient by the
equidistributional property of the nodal set of eigenfunctions. In
Section 3 we reduce by rescaling the gradient estimate from above
to that of an almost harmonic function. In Section 4, we prove the
gradient estimate from above for the almost harmonic function by
Yau's gradient estimate and the Green function for compact
manifolds with boundary. We also in this section sketch a proof of
the above estimates (\ref{equ:spegrad}) and (\ref{equ:2infgrad})
by the same argument with Theorem, which is different from those
in \cite{XuX} and \cite{XuB}.  In the last section we propose a
conjecture that a similar gradient estimate as Theorem holds for
a compact Riemannian manifold with boundary.\\

\section{Nodal sets of eigenfunctions}

\quad \ \ The nodal set of an eigenfunction $e_\l$ is the zero set
\[Z_{e_\l}:=\{x\in M: e_\l(x)=0\}.\]
A connected component of the open set $M\backslash Z_{e_\l}$ is
called a nodal domain of the eigenfunction $e_\l$. The following
fact, due to J. Br{\" u}ning \cite{Br}, is critical in deducing
the lower bound estimate. The reader can also find a proof in
English in Theorem 4.1 of S. Zelditch
\cite{Z}.\\

\nd {\bf Fact 1} {\it Under the assumption of Theorem, there
exists a constant $C$ such that each geodesic ball of radius
$C/\l$  in $M$ must
intersect the nodal set $Z_{e_\l}$ of $e_\l$.}\\

\nd {\sc Proof of the lower bound part of Theorem}\quad Take a
point $x$ in $M$ such that $|e_\l(x)|=\|e_\l\|_\infty$. Then by
Fact 1, there exists a point $y$ in the ball $B(x,\,C/\l)$ with
center $x$ and radius $C/\l$ such that $e_\l(y)=0$. We may assume
$\l$ so large that there exists a geodesic normal chart $(r,
\theta)\in [0,\,C/\l]\times {\Bbb S}^{n-1}(1)$ in the ball
$B(x,\,C/\l)$. By the mean value theorem, there exists a point $z$
on the geodesic segment connecting $x$ and $y$ such that
\[\left|\frac{\p e_\l}{\p r}(z)\right|\geq
\frac{\l}{C}|e_\l(x)|=\frac{\l}{C}\|e_\l\|_\infty.\]

\section{Eigenfunctions on the wavelength scale}
\quad \ \ In this section we review quickly the following principle:\\

\nd {\it On a small scale comparable to the wavelength $1/\l$, the
eigenfunction
$e_\l$ behaves like a harmonic function.}\\

\nd This principle was developed in H. Donnelly  and C. Fefferman
\cite{DF1} \cite{DF2} and N. S. Nadirashvili \cite{Na} and was
used extensively there. Our setting of the principle in the
following is borrowed from D. Mangoubi \cite{Ma1}, where he
applied it to studying the geometry of nodal domains of
eigenfunctions.

Fix an arbitrary point $p$ in $M$ and choose a geodesic normal
chart $(B(p,\d),\, x=(x_1,\cdots, x_n))$ with center $p$ and
radius $\d>0$ depending only on $M$. In this chart, we may
identify the ball $B(p,\d)$ with the $n$-dimensional Euclidean
ball ${\Bbb B}(\d)$ centered at the origin $0$, and think of the
eigenfunction $e_\l$ as a function in ${\Bbb B}(\d)$.  Our aim is
to show the inequality $|\nabla e_\l(p)|\leq C\l\|e\|_\infty$.
Actually we plan to prove a slightly stronger one
\begin{equation}
\label{equ:upper} \sum_j\,\left|\frac{\p e_\l}{\p
x_j}(0)\right|\leq C \l\|e_\l(x)\|_{L^\infty\left({\Bbb
B}(1/\l)\right)},
\end{equation}
which we remember that $g_{ij}(0)=\d_{ij}$ in the normal chart.
The rough idea is to consider the rescaled function
$\phi_\l(y)=e_\l(y/\l)$ in the ball ${\Bbb B}(1)$ instead of the
restriction of the eigenfunction $e_\l(x)$ to the ball ${\Bbb
B}(1/\l)$, where we assume $\l$ so large that $1/\l<\d/2$. The
above estimate is equivalent to its rescaled version
\begin{equation}
\label{equ:upperscale} \sum_j\,\left|\frac{\p \phi_\l}{\p
y_j}(0)\right|\leq C \|\phi_\l(y)\|_{L^\infty\left({\Bbb
B}(1)\right)}.
\end{equation}
On the other hand, rescaling the eigenfunction equation $\Delta_M
e_\l=\l^2 e_\l$ tells us that the function $\phi_\l$ behaves like
a harmonic function. The details are given in what follows.

The eigenfunction equation $\Delta_M e_\l=\l^2 e_\l$ in ${\Bbb
B}(1/\l)$ can be written as

\[-\frac{1}{\sqrt{g}}\sum_{i,j}\,\p_{x_i}\left(g^{ij}\sqrt{g}\p_{x_j}
e_\l\right)=\l^2 e_\l.\] Hence the function $\phi_\l$ satisfies
the rescaled equation in ${\Bbb B}(1)$,
\begin{equation}
\label{equ:resc}
-\frac{1}{\sqrt{g_\l}}\sum_{i,j}\,\p_{y_i}\left(g^{ij}_\l\sqrt{g_\l}\p_{y_j}
\phi_\l\right)=\phi_\l, \end{equation} where
$g_{ij,\,\l}(y)=g_{ij}(y/\l)$, $g^{ij}_\l(y)=g^{ij}(y/\l)$ and
$\sqrt{g_r}(y)=(\sqrt{g})(y/\l)$. We note that the above equation
coincides with
\[\D_\l \phi_\l=\phi_\l,\]
where $\Delta_\l$ is the positive Laplace-Beltrami operator with
respect to the metric $g_{ij,\l}$ on ${\Bbb B}(1)$. We observe
that the scaling $x\mapsto y=\l x$ gives an isometry from $({\Bbb
B}(1/\l),\,\l^2 g_{ij})$ onto $({\Bbb B}(1),\,g_{ij,\,\l})$.
Therefore, the open Riemannian manifold $({\Bbb
B}(1),\,g_{ij,\,\l})$ has uniformly bounded sectional curvature
for all $\l\geq 1$. This fact will be crucial in next section.


\section{Estimates from above}

\quad We prove (\ref{equ:upperscale}) in this section by using the
notations in Section 3.

{\it Step 1} \quad Recall that $\D_\l \phi_\l=\phi_\l$ in
$\BB(1)$. We can write the function $\phi_\l$ as the sum
$\phi_\l=u_\l+v_\l$ such that the functions $u_\l$ and $v_\l$
satisfy the following boundary-value problems:
\begin{equation}
\label{equ:bv} \left\{
\begin{array}{rcl}
\D_\l u_\l=0  & {\rm in}\
\BB(1), \\
u_\l=\phi_\l  & {\rm on}\ \p \BB(1),
\end{array}
\right. \quad\quad \left\{
\begin{array}{rl}
\D_\l v_\l=\phi_\l& {\rm in}\
\BB(1), \\
v_\l=0 &   {\rm on}\ \p \BB(1),
\end{array}
\right.
\end{equation}

{\it Step 2}\quad By the gradient estimate for harmonic functions
in page 21 of R. Schoen and S.-T. Yau \cite{SY} and the maximum
principle, there holds
\[\sum_j\,\left|\frac{\p u_\l}{\p
y_j}(0)\right|\leq C \|u_\l(y)\|_{L^\infty\left({\Bbb
B}(1)\right)}\leq C\|\phi_\l(y)\|_{L^\infty\left({\Bbb
B}(1)\right)}.\] Here we also use the fact that the Ricci
curvature tensor of the metrics $g_{ij,\l}$ are uniformly bounded
for all $\l\geq 1$.

{\it Step 3}\quad Let ${\Bbb G}(y,\,z)$ be the Green function on
the compact manifold $(\overline{\BB(1)},\,g_{ij,\l})$ with smooth
boundary such that in the sense of distribution $\D_{\l,\,z}{\Bbb
G}(y,\,z)$ equals the Dirac delta function $\d_y(z)$ at $y$ and
${\Bbb G}(y,\,z)$ vanishes whenver $y$ or $z$ belongs to the
boundary $\p \BB(1)$. Remember that the sectional curvatures of
$g_{ij,\l}$ are uniformly bounded for all $\l\geq 1$. Carefully
checking the explicit construction of the Green kernel in pages
106-113 in T. Aubin \cite{Au}, we find that
\[|(\nabla_y {\Bbb G})(0,\,z)|\leq C\,d(0,\,z)^{1-n}.\]
Taking gradient at $0$ in the Green formula for $v_\l$
\[v_\l(y)=\int_{\BB(1)}\, {\Bbb G}(y,\,z)\phi_\l(z)\,dV(z)\]
gives the estimate $|(\nabla\,v_\l)(0)|\leq
C\|\phi_\l\|_{L^\infty(\BB(1))}$. Remember that
$g_{ij,\l}(0)=\d_{ij}$ so that $|(\nabla\,v_\l)(0)|$ is comparable
to $\sum_j\,|(\p v_\l/\p y_j)(0)|$. The proof is completed.\\

Following the above arguments, we conclude this section by
sketching an alternative proof of the estimates
(\ref{equ:spegrad}) and (\ref{equ:2infgrad}), which were proved by
X. Xu \cite{XuX} and the second author \cite{XuB} by using the
maximum principle and the parametrix of the wave kernel,
respectively. We also use the notations above. Take a square
intergable function $f$ on $M$ and consider the smooth function
$u:=\chi_\l f$, and the Poisson equation $v:=\D_M\,u$ in the
geodesic normal chart $\BB(1/\l)$ centered at a given point $p$ in
$M$ and with radius $1/\l$. Rescaling the ball $\BB(1/\l)$ with
metric $g_{ij}$ to $\BB(1)$ with metric $g_{ij,\,\l}$, and using
the $L^\infty$ estimate \[\|\chi_\l f\|_{L^\infty}\leq
C\l^{(n-1)/2}\|f\|_2\] implied by the local Weyl law and the
Cauchy-Schwarz inequality, we reduce the question to showing
\begin{equation*}
\sum_{j}\,\left|\frac{\p u_\l}{\p y_j}(0)\right|\leq
C\|u\|_{L^\infty(\BB(1/\l))}+C\l^{-2}\,\|v_\l\|_{L^\infty(\BB(1))},
\end{equation*}
where $u_\l(y)=u(y/\l)$ and $v_\l(y)=v(y/\l)$. But the last
inequality follows from the arguments in Steps 1-3 above. The
proof is completed.

\section{Compact manifold with boundary}

D. Grieser \cite{Gr} and C. D. Sogge \cite{S3} obtained a similar
sup norm estimate as (\ref{equ:sup}) for $\chi_\l$ associated with
either the Dirichlet or the Neumann Laplacian on a compact
Riemannian manifold $N$ with boundary. X. Xu \cite{X} \cite{X2}
obtained similar gradient estimates as (\ref{equ:2infgrad}) for
both the Dirichlet and the Neumann boundary value problems on $N$
by using a clever maximum principle and the results of Grieser and
Sogge. However, we could not use the maximum principle argument in
their papers to deduce more refined gradient estimate for a single
eigenfunction
on $N$. We think that new ideas should be introduced to resolve the following conjecture. \\

\nd {\bf Conjecture} {\it Let $N$ be a compact Riemannian manifold
with smooth boundary and $e_\l(x)$ be an eigenfunction of either
the Dirichlet or the Neumann Laplacian $\Delta_N$ on $N$:
$\Delta_N\,e_\l=\l^2\,e_\l$. Then there exists a positive constant
$C$ depending only on $N$ such that for all $\l\geq 1$ there
holds}
\[\l\|e_\l\|_\infty/C\leq \|\nabla\,e_\l\|_\infty\leq
C\l\|e_\l\|_\infty.\]

We should point out that the lower bound estimate for the
Dirichlet eigenfunction $e_\l$,
\[\|\nabla\,e_\l\|_\infty\geq C\l \|e_\l\|_\infty,\]
follows easily from the similar argument in Section 2. We also
observe that the upper bound estimate $|\nabla\,e_\l(x)|\leq C\l
\|e_\l\|_\infty$ for every $x$ outside the boundary layer $\{x\in
N:d(x,\,\p N)\leq C/\l\}$ follows from the argument in Sections 3
and 4. The authors will discuss  this question in a future
paper.\\

\nd \begin{center} {\sc Acknowledgements}
\end{center}

The authors would like to thank Professor Qing Chen, Professor Yu
Safarov, Professor C. D. Sogge and Professor John Toth for their
valuable conversations and comments. The second author would like
to express his sincere gratitude to Professor Qing Chen and
Professor Xiuxiong Chen for their constant encouragement and
forceful prodding.

\end{document}